\documentclass{article}

\usepackage{amsmath,amsfonts,amsthm,amstext}
\usepackage{amsmath,amssymb,indentfirst,epsfig,multicol}
\usepackage{amssymb,latexsym}

\usepackage{epsfig}
\usepackage{graphicx}
\newcommand{\Z}{\mathbb Z}
\newcommand{\R}{\mathbb R}
\newcommand{\cqd}{\hfill $\rule{2.5mm}{2.5mm}$}
\newcommand{\ind}{\mathop{\mathrm{ind}}}
\newcommand{\cohomdim}{\mathop{\mathrm{cohom.dim}}}
\title{$(H,G)$-coincidence theorems for manifolds and a  topological Tverberg type theorem for any natural number $r$}

\author{Denise de Mattos\thanks{The author was supported in part by FAPESP of Brazil
Grant numbers 12/24454-8 and 2013/24845-0},  
 Edivaldo L. dos Santos\thanks{The author was supported in part by FAPESP of Brazil
Grant numbers 12/24454-8 } and Taciana O. Souza\thanks {The author was supported by CAPES and
CNPq} }

\date{}

\newtheorem{theorem}{Theorem}[section]
\newtheorem{rema}[theorem]{Remark}
\newtheorem{lema}[theorem]{Lemma}

\newtheorem{defi}[theorem]{Definition}

\newtheorem{cor}[theorem]{Corollary}

\newcommand{\dist}{{\rm{dist}}}
\begin{document}
\maketitle

\begin{abstract} Let $X$ be a paracompact space, let $G$ be a finite group
acting freely on $X$ and let $H$ a cyclic subgroup of $G$ of  prime
order $p$. Let $f:X\rightarrow M$ be a continuous map where $M$ is a
connected $m$-manifold (orientable if $p>2$) and $f^* (V_k) = 0$,
for $k\geq 1$, where $V_k$ are the $Wu$ classes of $M$. Suppose that
$\ind X\geq n> (|G|-r)m$, where $r=\frac{|G|}{p}$. In this work, we 
estimate the cohomological dimension of the set $A(f,H,G)$ of
$(H,G)$-coincidence points of $f$. Also,  we  estimate the index of a $(H, G)$-coincidence set in the case that $H$ is a $p$-torus subgroup of a particular group $G$ and as application we prove a topological Tverberg type theorem for any natural number $r$. Such result is a weak version of the famous topological Tverberg conjecture, which  was proved recently, fail for all $r$ that are not prime powers. Moreover, we obtain a  generalized Van Kampen-Flores type  theorem for any  natural number $r$.

\vspace{0.1cm}

\end{abstract}

\begin{section}{Introduction}

Let $G$ be a finite group which acts freely on a space $X$ and let
$f:X\rightarrow Y$ be a continuous map from $X$ into another space
$Y$. If $H$ is a subgroup of $G$, then $H$ acts on the right on each
orbit $Gx$ of $G$ as follows: if $y\in Gx$ and $y=gx$, with $g\in
G$, then $h\cdot y=gh^{-1}x$.  A point $x\in X$ is said to be a
$(H,G)$- \emph{coincidence point of $f$} (as introduced by Gon\c
calves, Pergher and Jaworowski in \cite{Pergher}) if $f$ sends every
orbit of the action of $H$ on the $G$-orbit of $x$ to a single
point. Of course, if $H$ is the trivial subgroup, then
every point of $X$ is a $(H,G)$-coincidence. If $H=G$, this is the
usual definition of $G$-coincidence, that is, $f(x)=f(gx)$, for all
$g\in G$. If $G=\mathbb{Z}_p$ with p prime, then a nontrivial
$(H,G)$-coincidence point is a $G$-coincidence point. Let us denote
by $A(f,H,G)$ the set of all $(H,G)$-coincidence points. Borsuk-Ulam
type theorems consists in estimating the cohomological dimension of
the set $A(f,H,G)$. Two main directions considered of this problem
are either when the target space $Y$ is a manifold or $Y$ is a CW
complex. In the first direction are the papers of Borsuk
\cite{Borsuk} ( the classical theorem of Borsuk-Ulam, for $H=G = \mathbb{Z}_2$,
$X = S^ {n}$ and $Y = R^{n}$), Conner and Floyd \cite{Conner}
(for $H=G = \mathbb{Z}_2$ and $Y$ a $n$-manifold), Munkholm \cite{Munkholm2}
(for $H=G = \mathbb{Z}_p$, $X = S^ {n}$ and $Y=R^{m}$), Nakaoka \cite{Nakaoka}
(for $H=G= \mathbb{Z}_p$, $X$ under certain (co)homological conditions and
$Y$ a $m$-manifold) and the following more general version proved by
Volovikov \cite{Volovikov} using the index of a free $\mathbb{Z}_p$-space
$X$ ($\ind X$, see Definition \ref{index} ):

\vspace{0.2cm}

\noindent {\bf Theorem A.}\cite[Theorem 1.2]{Volovikov} {\it Let $X$
be a paracompact free $\mathbb{Z}_p$-space of  $\ind X \geq n$, and $f: X
\to M$ a continuous mapping  of $X$ into an $m$-dimensional
connected manifold $M$ (orientable if $p
> 2$). Assume that:

\vspace{0.2cm}

\noindent {\rm (1)} $f^*(V_{i}) = 0$ for $i \geq 1$, where the $V_i$
are the Wu classes of $M$; and

\vspace{0.2cm}

\noindent {\rm (2)} $n>m(p- 1)$.

\vspace{0.2cm}

Then the $\ind A(f)\geq n - m(p - 1)> 0$.}

\vspace{0.2cm}

In the second direction are the papers of Izydorek and Jaworowski
\cite{Izydorek1} (for $H=G= \mathbb{Z}_p$, $X=S^{n}$ and $Y$ a CW-complex ),
Gon\c calves and Pergher \cite{Pergher2} (for $H=G= \mathbb{Z}_p$, $X=S^{n}$
and $Y$ a CW-complex ) and for proper nontrivial subgroup $H$ of
$G$, Gon\c calves, Jaworowski and Pergher \cite{Pergher} (for
$H= \mathbb{Z}_p$ subgroup of a finite group $G$, $X$ an homotopy sphere and
$Y$ a CW-complex) and   Gon\c calves, Jaworowski, Pergher and
Volovikov \cite{Daciberg}(for $H= \mathbb{Z}_p$ subgroup of a finite group
$G$, $X$ under certain (co)homological assumptions and $Y$ a
CW-complex).

In this work, considering the target space $Y=M$ a manifold and $H$
a proper nontrivial subgroup of $G$, we prove the following
formulation of the Borsuk-Ulam theorem for manifolds in terms of
$(H,G)$-coincidence.

\begin{theorem}\label{teo1} Let $X$ be a paracompact space of $ind X\geq n$ and let $G$ be a finite group
acting freely on $X$ and $H$ a cyclic subgroup of $G$ of  prime
order $p$. Let $f:X\rightarrow M$ be a continuous map where $M$ is a
connected $m$-manifold (orientable if $p>2$) and $f^* (V_k) = 0$,
for all $k\geq 1$, where $V_k$ are the $Wu$ classes of $M$. Suppose
that  $\ind X\geq n> (|G|-r)m$ where $r=\frac{|G|}{p}$, then $\ind
A(f,H,G)\geq n-(|G|-r)m$. Consequently, $$\cohomdim A(f,H,G)\geq
n-(|G|-r)m>0.$$
\end{theorem}

\vspace{0.5cm}

Let us observe that if $H=G= \mathbb{Z}_p$, we have
$(|G|-r)m=(p-1)m$ and therefore Theorem \ref{teo1} generalizes
Theorem A. obtained by Volovikov in \cite{Volovikov} in terms of
$(H,G)$-coincidence. For the case $n=(|G|-r)m$, $p$ an odd prime, if
we consider $X$ a $mod p$ homology $n$-sphere in the Theorem
\ref{teo1}(in this case,  the continuous  map $f$ can be arbitrary),
it follows from \cite[Theorem 8]{Nakaoka} that $A(f,H,G)\neq
\emptyset$. Also, Theorem \ref{teo1} is a version for manifolds of
the main result due to Gon\c{c}alves, Jaworowski and Pergher in
\cite{Pergher}.

\vspace{0.2cm}

 Also, we prove the following  nonsymmetric theorem for
$(H,G)$-coincidences which is a version for manifolds of the main
theorem in \cite{DM2}.

\begin{theorem}\label{teo2}Let $X$ be a compact Hausdorff space, let $G$ be a finite group
acting freely on $S^n$ and let $H$ be a cyclic subgroup of $G$ of
order prime $p$. Let $\varphi : X\rightarrow S^n$ be an essential
map \footnote{A map $\varphi : X\rightarrow S^n$ is said to be an
essential map if $\varphi$ induces nonzero homomorphism $\varphi^* :
H^n (S^n; {\mathbb Z}_p)\rightarrow  H^n (X; {\mathbb Z}_p)$.} and
let $f:X\rightarrow M$ be a continuous map where $M$ is a connected
$m$-manifold (orientable if $p>2$) and $f^* (V_k) = 0$, for all
$k\geq 1$, where $V_k$ are the $Wu$ classes of $M$. Suppose that $n>
(|G|-r)m$, then
$$\cohomdim A_{\varphi}(f,H,G)\geq n-(|G|-r)m,$$
where $r=\frac{|G|}{p}$ and $A_{\varphi} (f,H,G)$ denotes the
$(H,G)$-coincidence points of $f$ relative to an essential map
$\varphi:X\to S^{n}$.
\end{theorem}

In Section 5,  we give similar estimate  in the case that  $H$ is a $p$-torus subgroup of a particular group $G$ and as application, we prove a topological Tverberg type theorem for any natural number, which is a weak version of the famous topological Tverberg conjecture. Moreover, we obtain a  generalized Van Kampen-Flores type  theorem for any interger $r$.

\end{section}

\section{Preliminaries}
We start by introducing some definitions as follows.

\begin{subsection}{The $\mathbb {Z}_p$-index}

We suppose that the cyclic group $\mathbb Z_p$ acts freely on a
paracompact Hausdorff space $X$, where $p$ is a prime number and we
denote by $[X]^*$ the space of orbits of $X$ by the action of
$\mathbb Z_p$. Then, $X\rightarrow [X]^*$ is a principal $\mathbb
Z_p$-bundle and we can consider a classifying  map $c:[X]^*
\rightarrow B \mathbb Z_p$.

\begin{rema} \rm It is well known that if $\hat{c}$ is another classifying map
for the principal $\mathbb{Z}_{p}$-bundle $X\to X^{*}$, then there
is a homotopy between $c$ and $\hat{c}$.
\end{rema}

\begin{defi}\label{index} We say that the $\mathbb Z_p$-index of $X$ is greater
than or equal to $l$ if the homomorphism $$c^* : H^l (B \mathbb Z_p
; \mathbb Z_p) \rightarrow H^l ([X]^* ; \mathbb Z_p)$$ is
nontrivial. We say that the $\mathbb Z_p$-index of $X$ is equal to
$l$ if it is greater or equal than $l$ and, furthermore, $c^* : H^i
(B \mathbb Z_p ; \mathbb Z_p) \rightarrow H^i ([X]^* ; \mathbb Z_p)$
is zero, for all $i\geq l+1$.

We denote the $\mathbb Z_p$-index of $X$ by $\ind X$.
\end{defi}
\end{subsection}

%%%%%%%%%%%%%%%%%%%%%%%%%%%%%%%%%%%%%%%%%%%%%%%%%%%%%%%%%%%%%%%%%%%%%%%%%%%%%%%%%5
\subsection{Wu classes}

The \textit{total Wu class} of a manifold $M$ is defined as the formal sum
$$v(M) = 1 + v_1 (M) + v_2 (M) + \cdots + v_k (M) + \cdots$$
where $v_k (M)$  is the $k$-th Wu class of $M$, $k = 1, 2, \ldots$ (see \cite{Milnor2}). Let $p>2$ be a prime, using the \textit{total reduced power} $$P = P^0 + P^1 + P^2 + \cdots + P^k +  \cdots$$ and the equation $$ \langle v_k (M) \smile x , [M] \rangle = \langle P^k (x), [M] \rangle $$ we obtain the formula $$\langle v(M) \smile u, [M] \rangle = \langle P (u), [M] \rangle $$ for all $u \in H^*_c (M; \mathbb{Z}_p).$ For $p=2$ we have a similar formula $$\langle v(M) \smile u, [M] \rangle = \langle Sq (u), [M] \rangle$$ for all $u \in H^*_c (M; \mathbb{Z}_2)$, where $$Sq = Sq^0 + Sq^1 + Sq^2 \cdots + Sq^k + \cdots$$ is the \textit{total Steenrod square}.

\section{The $Wu$ class for product of
manifolds} 
Let $W$ and $M$ be manifolds with
dimensions $w$ and $m$ respectively, $w\geq m$, both orientables
if $p>2$. In the next lemma, we obtain a characterization of the $Wu$
class of the product of $W$ and $M$.

\begin{lema} \label{lema0} Let $W$ and $M$ be a connected manifolds, both orientables
if $p>2$. Then, the total $Wu$ class of  $W\times M$, is given
by:
\begin{eqnarray} v(W)\otimes v(M)
\end{eqnarray}
where $v(W)$ and $v(M)$ are the total $Wu$ classes of $W$ and $M$ respectively.

\end{lema}

\noindent {\it Proof.} Let $p>2$ be a prime number. Let $z = w\otimes u$ an element of $ H^*_c (W \times M; \mathbb{Z}_p)$ then

$$
\begin{array}{lll}

\langle v(W)\otimes v(M) \smile z , [W \times M] \rangle & = & \langle v(W) \smile w \otimes v(M) \smile u , [W \times M] \rangle \\
& = & \langle P(w) \otimes P(u) , [W \times M] \rangle \\
& = & \langle P(w \otimes u) , [W \times M] \rangle \\
& = & \langle P(z) , [W \times M] \rangle \\
& = & \langle v(W \times M) \smile z , [W \times M] \rangle \\
\end{array}
$$
Therefore by uniqueness of the Wu class we conclude that the total Wu class of $W \times M$ is given by $v(W \times M) = v(W) \otimes v(M)$. By a similar argument the total Wu classes are obtained for $p=2$, in this case are used the total Steenrod square. \cqd

\vspace{0.3cm}

\section{$(H,G)$-Coincidence theorems for manifolds}

Now, we denote by $a_1,\ldots,a_r$ a set of representatives of the
left lateral classes of $G/H$. We define, for each $i=1,\ldots,r$,
$g_i : X\rightarrow X$ by $g_i (x)=a_i x$. Consider the map $F:
X\rightarrow M^r$ defined by $$F = (f_1 \times \ldots \times
f_r )\circ D ,$$ where $D : X\rightarrow X^r $ is the diagonal map
and $f_i = f\circ g_i$. We prove the following

\begin{lema}\label{lema1} If $f^* (v_k (M)) = 0$, for all $k\geq
1$, where  $v_k (M)$ are the $Wu$ classes of $M$, then $F^* (v_k (M^r))=0$, for all $k\geq 1$, where  $v_k (M^r)$ are the $Wu$ classes of $M^r$.
\end{lema}

\noindent {\it Proof.} It suffices to show that $(f_1 \times \ldots \times
f_r )^* (v_k (M^r)) = 0$, for $k\geq 1$. If $r=1$, then $F=f_1$ and $f_1 ^* (v_k (M)) = g_1 ^* \circ f^* (v_k (M)) = 0$. 

Let us denote by $$p_1 : M^{r-1} \times M \to  M^{r-1}, \, p_2 :M^{r-1} \times M \to  M $$ 
$$q_1 : X^{r-1} \times X \to  X^{r-1},  \, q_2 : X^{r-1} \times X \to  X$$
the natural projections. If $r\geq2$, we have $$(f_1 \times \ldots \times f_{r-1}) \circ q_1  = p_1 \circ (f_1 \times \ldots \times f_{r})$$ 
$$f_r \circ q_2 = p_2 \circ (f_1 \times \ldots \times f_{r}).$$ 
Since, by Lemma \ref{lema0}, $\displaystyle v_k (M^{r-1} \times M) = \sum_{s=0}^{k} v_{s} (M^{r-1})\times v_{k-s} (M)$ 
and assuming inductively that $(f_1 \times \ldots \times f_{r-1})^* (v_{s} (M^{r-1})) = 0$, for $s\geq 1$, 
we conclude that

$$
\begin{array}{l}
 \ \ (f_1 \times \ldots \times f_{r})^* (v_k (M^{r-1} \times M)) = \\
= (f_1 \times \ldots \times f_{r})^* \left(  \displaystyle \sum_{s=0}^{k} v_{s} (M^{r-1}) \times v_{k-s} (M)\right) \\
=  \displaystyle \sum_{s=0}^{k} (f_1 \times \ldots \times f_{r})^* (p^*_1 (v_{s} (M^{r-1}))) \smile (f_1 \times \ldots \times f_{r})^* (p_2^* (v_{k-s} (M))) \\
= \displaystyle \sum_{s=0}^{k} q_1^{*} \circ (f_1 \times \ldots \times f_{r-1})^*(v_{s} (M^{r-1})) \smile q_2^{*} \circ g_r^{*} \circ f^* (v_{k-s}(M)) \\
 = 0. \cqd
\end{array} 
$$

\subsection{Proof of Theorem \ref{teo1} and its consequences}

\noindent In this section, we present the proofs of Theorems
\ref{teo1} and its consequences, as follows.

\vspace{0.2cm}

\begin{Proof}{\it{\, of Theorem \ref{teo1}}}. We consider the map $F$ defined previously. We have
$$A(f,H,G)\supset A_F=\{x\in X : F(x) = F(hx), \forall h\in H \}.$$

In fact, let $x$ be a point in the set $A_F$, then
$$(f(a_1 x),\ldots,f(a_r x))= (f(a_1 h x),\ldots, f(a_r h x)),$$
for all $h \in H$. Thus, $f(a_i x)= f(a_i hx)$, for all $h\in H$ and
$i=1,\ldots,r$. According to the definition of the action of $H$ on
the orbit $Gx$, $h^{-1}\cdot a_i x := a_i (h^{-1})^{-1} x=a_{i}hx\in
a_{i}Hx$, for $i=1,\ldots, r$. Thus, $f$ collapses each orbit
$a_{i}Hx$ determined by the action of $H$ on $a_i x$, for
$i=1,\ldots, r$, therefore $x\in A(f,H,G)$.

Now we observe that $H\cong \mathbb Z_p$ acts freely on $X$ by
restriction and by hy\-po\-thesis $\ind X\geq n> n-(p-1)rm$. By
Lemma \ref{lema1},  $F^* (v_k) = 0$, for all $k\geq 1$, where $v_k$
are the $Wu$ classes of $M^r$. Thus, according to \cite[Theorem
$1.2$]{Volovikov}
$$\ind A_F\geq  n-(p-1)rm = n-(|G|-r)m.$$

Let us consider the inclusion $i: A_F\rightarrow A(f,H,G)$, which is
an equivariant map, and so it induces $\overline{i} : [A_F]^*
\rightarrow [A(f,H,G)]^*$ a map between the orbit spaces. Therefore,
if $c : [A(f,H,G)]^* \rightarrow B \mathbb Z_p$ is any classifying
map, we have that $c \circ \overline{i} : [A_F]^* \rightarrow B
\mathbb Z_p$ is a classifying map. Thus,
$$\ind A(f,H,G)\geq \ind A_F \geq n- (|G|-r)m,$$
  \cqd
\end{Proof}

\begin{cor} \label{cor 1}Let $X$ be a paracompact space and let $G$ be a finite group
acting freely on $X$. Let $M$ be a orientable $m$-manifold, and $p$
a prime number that divide $|G|$. Suppose that $ind X \geq n>
(|G|-r)m$, where $r= \frac{|G|}{p}$. Then, for a continuous map $f:X
\rightarrow M$ such that $f^* (V_k) = 0$, for all $k\geq 1$, where
$V_k$ are the $Wu$ classes of $M$, there exists a non-trivial
subgroup $H$ of $G$, such that $$\cohomdim A(f,H,G)\geq
n-(|G|-r)m.$$
\end{cor}

\begin{Proofl} Let $p$ be a prime number such that divide $|G|$. By
Cauchy Theorem, there is a cyclic of order $p$ subgroup $H$ of $G$ .
Then, we apply Theorem \ref{teo1}.

\cqd
\end{Proofl}

%%%%%%%%%%%%%%%%%%%%%%%%%%%%%%%%%%%%%%%%%%%%%%%%%%%%%%%%%%%%%%%%%%%%%%%%%%%%%%

\vspace{0.2cm}

\begin{rema}\rm Let us observe that, if $f^{*} : H^i (M; \mathbb{Z}_p)\rightarrow H^i (X; \mathbb{Z}_p)$ is trivial, for $i\geq 1$, and
$p$ is the smallest prime number dividing $|G|$, then $r=
\frac{|G|}{p}\geq \frac{|G|}{q}$, where $q$ can be any other prime
number dividing $|G|$. Thus, $n> (|G|- \frac{|G|}{q})m$, therefore
for each prime number $q$ dividing $|G|$, there exists a cyclic
subgroup of order $q$, $H_q$ of $G$ such that $\ind A(f,H_q ,G)\geq
n-(|G|-r)m$.
\end{rema}

%%%%%%%%%%%%%%%%%%%%%%%%%%%%%%%%%%%%%%%%%%%%%%%%%%%%%%%%%%%%%%%%%%%%%%%%%%%%%%%%%%%%%%%%%%%%%%%%%

\vspace{0.2cm}

The following theorem is a version for manifolds of the main result in \cite{Pergher}

\begin{theorem} Let $G$ be a finite group which acts freely on $n$-sphere $S^n$
and let $H$ be a cyclic subgroup of $G$  of  prime order $p$. Let
$f:S^{n}\to M$ be a continous map where $M$ be a $m$-manifold
(orientable if $p>2$). If $n>(|G|-r)m$ where $r=\frac{|G|}{p}$, then
$$\cohomdim (A(f, H, G))\geq n-(|G|-r)m.$$
\end{theorem}

\begin{Proofl}  Since $n>(|G|-r)m\geq m$, $f^{*}(V_{k})=0$, for all $k\geq 1.$
Moreover, $\ind S^{n}=n$ and thus we apply the Theorem \ref{teo1}.
\cqd
\end{Proofl}

\subsection{Proof of Theorem \ref{teo2}}

Now, let us consider $X$ a compact Hausdorff space and an essential
map $\varphi : X\rightarrow S^n$. Suppose $G$ be a finite group de
order $s$ which acts freely on $S^n$ and $H$ be a subgroup of order
$p$ of $G$. Let $G = \{g_1, ..., g_s \}$ be a fixed enumeration of
elements of $G$, where $g_1$ is the identity of $G$. A nonempty
space $X_{\varphi}$ can be associated with the essential map
$\varphi : X\rightarrow S^n$ as follows:

$$X_{\varphi}= \{(x_1,..., x_s)\in X^s : g_i \varphi (x_1) = \varphi (x_i), i = 1,..., s \},$$
where $X^s$ denotes the $s$-fold cartesian product of $X$. The set
$X_{\varphi}$ is a closed subset of $X^s$ and so it is compact. We
define a $G$-action on $X_{\varphi}$ as follows: for each $g_i \in
G$ and for each $(x_1,..., x_s) \in X_{\varphi}$,
$$g_i (x_1,..., x_s) = (x_{\sigma_{g_i} (1)},..., x_{\sigma_{g_i} (s)}),$$
where the permutation $\sigma_{g_i}$, is defined by $\sigma_{g_i}
(k) = j$, $g_{k} g_{i} = g_{j}$. We observe that if $x = (x_1, ...,
x_s)\in
 X_{\varphi}$ then $x_i\neq x_j$, for any $i\neq j$ and therefore
$G$ acts freely on $X_{\varphi}$.

Let us consider a continuous map $f: X\rightarrow M$, where $M$ is a
topological space and $\widetilde{f} : X_{\varphi}\rightarrow M$
given by $\widetilde{f} (x_1,..., x_s) = f(x_1)$,

\begin{defi} The set $A_{\varphi} (f,H,G)$ of $(H,G)$-coincidence
points of $f$ relative to $\varphi$ is defined by
$$A_{\varphi} (f,H,G) = A(\widetilde{f},H,G).$$
\end{defi}

\begin{Proof} {\it{\, of Theorem \ref{teo2}.}}
Let $\widetilde{f}: X_{\varphi}\rightarrow M$ given by
$\widetilde{f} (x_1, ..., x_r) = f(x_1)$, that is, $\widetilde{f} =
f\circ \pi_1$, where $\pi_1$ is the natural projection on the $1$-th
coordinate. By hypothesis,  $f^* (V_k) = 0$, for all $k\geq 1$,
where $V_k$ are the $Wu$ classes of $M$, then we have
$\widetilde{f}^* (V_k) = 0$, for all $k\geq 1$. Moreover, the
$\mathbb Z_p$-index of $X_{\varphi}$ is equal to $n$ by \cite{DM2}
Theorem $3.1$. In this way, $X_{\varphi}$ and $\widetilde{f}$
satisfy the hypothesis of Theorem \ref{teo1} which implies that the
$\mathbb Z_p$-index of the set $A(\widetilde{f}, H,G)$ is greater
than or equal to $n-(|G|-r)m$. By definition, $A_{\varphi} (f, H,G)=
A(\widetilde{f}, H,G)$, and then
$$\cohomdim A_{\varphi} (f, H,G)\geq n-(|G|-r)m.$$  \cqd
\end{Proof}

\vspace {0.2cm}

By a similar argument to that used in the proof of Corollary
\ref{cor 1} we have the following corollary of Theorem \ref{teo2}

\begin{cor}Let $X$ be a compact Hausdorff space and let $G$ be a finite group
acting freely on $S^n$. Let $M$ be a orientable $m$-manifold and $p$
a prime number dividing $|G|$. Suppose that $n> (|G|-r)m$, where $r=
\frac{|G|}{p}$. Then, for a continuous map $f:X \rightarrow M$, with
$f^* (V_k) = 0$, for all $k\geq 1$, where $V_k$ are the $Wu$ classes
of $M$, there exists a non-trivial subgroup $H$ of $G$, such that
$$\cohomdim A_{\varphi} (f, H,G)\geq n-(|G|-r)m.$$
\end{cor}

\section{Topological Tverberg type theorem }
The history of Tverberg theorem begins with a Birch's paper  (see \cite{Birch}) which contained the following conjecture

\begin{center}
{\it ``Any $(r - 1)(d+1) + 1$ points in $\mathbb{R}^d$ can be partitioned in $N$ subsets whose convex hulls have a common point ". }
\end{center}

The Birch's conjecture was proved by Helge Tverberg (see \cite{Tverberg}) and since then is known as Tverberg theorem. 

We note that the convex hull of $l + 1$ points in  $\mathbb{R}^d$ is the image of the linear map $\Delta_l \to \mathbb{R}^d$ that maps the $l + 1$ vertices of $\Delta_l$ to these $l + 1$ points. Thus the Tverberg theorem can be reformulated as follows:

\vspace{0.2cm}

\noindent {\bf Tverberg Theorem.} {\it Let $f$ be a linear map from the $N$-dimensional simplex $\Delta_N$ to $\mathbb{R}^d$. If $N = ( d + 1 )(r - 1)$ then there are $r$ disjoint faces of $\Delta_N$ whose images have a   common point.}

\vspace{0.2cm}

The following conjecture is a generalization of Tverberg Theorem to arbitrary continuous maps.

\vspace{0.2cm}

\noindent {\bf The topological Tverberg conjecture.} {\it Let $f$ be a continuous map from the $N$-dimensional simplex $\Delta_N$ to $\mathbb{R}^d$. If $N = ( d + 1 )(r - 1)$ then there are $r$ disjoint faces of $\Delta_N$ whose images have a  common point.} 

\vspace{0.2cm}

The topological Tverberg conjecture  was considered a central unsolved problem of topological combinatorics. For   a prime number $r$ the conjecture was proved by B\'ar\'any, Shlosman and  Sz\H{u}cs (\cite{Barany}) and it was extended for  a prime power $r$ by \"Ozaydin (unpublished) (\cite{Ozaydin}) and  Volovikov (\cite{Volovikov2}). This result is known as the {\it topological Tverberg thereom}. Recently, in \cite{Frick}, Frick presents surprising counterexamples  to the topological Tverberg conjecture for any $r$ that is not a power of a prime and dimensions $d\geq 3r +1$ (see also \cite{BFZ}).  Although, the conjecture is not  true for  an integer $r\geq 6$ that is not a prime power,  it is possible to prove a  weak version of the topological Tverberg conjecture, more precisely, in this paper we show that if $r$ is a natural number with prime factorization $r = p_1^{n_1}\cdots p_k^{n_k}$ then there is, for each $j = 1, \ldots, k$, a set with $r$ closed sides mutually disjoint of $\Delta_N$ which can be divided into $ \displaystyle \frac{r}{p_j^{n_j}}$ subsets, each one having  $p_j^{n_j}$ elements,  whose images have a common point. Specifically, we prove the following  {\it Topological Tverberg type theorem for manifolds and for any natural number $r$.}

\begin{theorem}\label{TverbergThm}
Let $d\geq 1$ a natural number. Consider a natural number $r$ with prime factorization $r = p_1^{n_1}\cdots p_k^{n_k}$ and set $N = (r-1)(d+1)$. Let $f: \partial {\Delta}_N \to M$ be a continuous mapping into a compact $d$-dimensional topological manifold. Then, for each $j = 1, \ldots, k$, among the sides of ${\Delta}_N$ there are $r=q_{j}r_j$, where $r_j = p_j^{n_j}$, and $q_j = \displaystyle \frac{r}{r_j}$, mutually disjoint closed sides  ${\sigma}_{1_1}$, ..., ${\sigma}_{1_{r_j}}$; $\ldots$ ;${\sigma}_{i_1}$, ..., ${\sigma}_{i_{r_j}}$; $\ldots$; ${\sigma}_{{q_{j}}_1}$, ..., ${\sigma}_{{q_{j}}_{r_{j}}}$, such that 
$$f({\sigma}_{i_1})\cap \cdots \cap f({\sigma}_{i_{r_j}}) \neq \emptyset, \,\, \mbox{for each}\,\, i=1,\ldots, q_{j}.$$
\end{theorem}

%Let $r$ be a integer with prime factorization $r = p_1^{n_1}\cdots p_k^{n_k}$. 

\begin{defi}[Index]
Let $p$ be a prime. We suppose the $p$-torus \linebreak $G = \mathbb{Z}_p ^k =  \mathbb{Z}_p \times \cdots \times  \mathbb{Z}_p$ (k factors) acting freely on a paracompact space $X$. The covering $X \to X/G$ is induced from the universal covering $EG \to BG$ by means of a classifying map $c: X/G \to BG$, defined uniquely up to homotopy. We say that the {\bf index} of $X$ is greater than or equal to $N$ (abbreviated by $\ind X \geq N$) if $c^*: H^N (BG;  \mathbb{Z}_p) \to H^N (X/G;  \mathbb{Z}_p)$ is a monomorphism.
\end{defi}

Consider $G = \mathbb{Z}_{p_1}^{n_1} \times \ldots \times \mathbb{Z}_{p_k}^{n_k}$, where $\Z_{p_j}^{n_j} = \mathbb{Z}_{p_j} \times \ldots \times \mathbb{Z}_{p_j}$ ($n_j$ factors), \linebreak 
$j = 1, \ldots, k$. We suppose $G$ acts freely on a paracompact space $X$. \linebreak Let $f: X \to M$ be a continuous mapping into a $d$-dimensional topological manifold.

\begin{lema}\label{lemaprime}
Let $M$ be a compact $d$-dimensional topological manifold orientable. Suppose the homomorphism $f^* : H^i (M; \mathbb{Z}_{p_j} ) \to H^i (X; \mathbb{Z}_{p_j})$ is trivial for $i\geq 1$,  and $\ind \, X \geq N \geq d(r - q_j)$, for each $j = 1, \ldots, k$, where $q_j = \displaystyle \frac{r}{p_j^{n_j}}$. Then   $$\ind  A\left(f, \mathbb{Z}_{p_j}^{n_j} , G \right) \geq N - d(q - q_j).$$
\end{lema}

\noindent {\it Proof.}  We denote by $a_1,\ldots,a_{q_j}$ a set of representatives of the
left lateral classes of $G/\mathbb{Z}_{p_j}^{n_j}$. We define, for each $i=1,\ldots,q_j$,
$g_i : X\rightarrow X$ by $g_i (x)=a_i x$. Consider the map $F:
X\rightarrow M^{q_j}$ defined by $$F = (f_1 \times \ldots \times
f_{q_j} )\circ D,$$  where $ (f_1 \times \ldots \times
f_{q_j} )(x_1, \ldots, x_{q_j}) = (f_1 (x_1), \ldots, 
f_{q_j} (x_{q_j})) $, $D : X\rightarrow X^{q_j} $ is the diagonal map
and $f_i = f\circ g_i$.

We have $F^*: H^i (M^{q_j}; \mathbb{Z}_{p_j} ) \to H^i (X; \mathbb{Z}_{p_j})$ trivial for $i\geq 1$, therefore the index of $A(F) = \{x \in X : F(x)= F(gx)  \,  \forall g \in \mathbb{Z}_{p_j}^{n_j}\}$ is greater  than or equal to $N - q_j d \left(p_j^{n_j} - 1 \right)$ (see \cite[Theorem 1]{Volovikov3}). Since $A(F)\subset A\left(f, \mathbb{Z}_{p_j}^{n_j} , G \right)$ and the inclusion $A(F)\hookrightarrow A\left(f, \mathbb{Z}_{p_j}^{n_j} , G \right) $ is an equivariant map we have $\ind A\left(f, \mathbb{Z}_{p_j}^{n_j} , G \right) \geq  A(F)$. Then
$$ \ind  A\left(f, \mathbb{Z}_{p_j}^{n_j} , G \right) \geq N - d(r - q_j).$$ \cqd

\vspace{0.3cm}

\noindent {\it Proof of Theorem \ref{TverbergThm}.} We consider the $CW$-complex $Y_{N,r}$ that consists of points $(y_1, \ldots, y_r)$, $y_i$ in the boundary $\partial \Delta_N$ of the simplex $\Delta_N$, that have mutually disjoint closed faces. It is known that for all natural numbers $r$ and $N$, where $N \geq r + 1$, $Y_{N,r}$ is $(N-r)$-connected (see \cite{Barany}). Let $G = \{g_1, \ldots, g_r\}$ be a fixed enumeration of elements of $G$, we define a $G$-action on $Y_{N,r} \subset ({\Delta}_N)^q$ as follows: for each $g_i \in G$ and for each $(y_1, \ldots, y_r) \in Y_{N,r}$
$$g_i (y_1,..., y_r) = (y_{\phi_{g_i} (1)},..., y_{\phi_{g_i} (r)}),$$
where the permutation $\phi_{g_i}$, is defined by $\phi_{g_i}
(k) = j$, $g_{k} g_{i} = g_{j}$. Then $G$ acts freely on $Y_{N,r}$, since $Y_{N,r}$ consists of points $(y_1, \ldots, y_r)$, $y_i \in \partial {\Delta}_N$ that have mutually disjoint closed faces.

Let $\tilde{f} : Y_{N,r} \to M$  given by $\tilde{f}(y_1, \ldots, y_r) = f(y_1)$, that is, $\tilde{f} = f \circ \pi_1$ where $\pi_1 : Y_{N,r} \to \partial {\Delta}^N$ is the projection on the 1-th coordinate.  Since \linebreak $\tilde{f}^* : H^i (M; \mathbb{Z}_{p_j}) \to H^i (Y_{N,r}; \mathbb{Z}_{p_j})$ is trivial for $i\geq 1$ and \linebreak $\ind  Y_{N,r} \geq d(r-1) > d(r - q_j)$ then, according to Lemma \ref{lemaprime}, the set $ A\left(\tilde{f}, \mathbb{Z}_{p_j}^{n_j} , G \right)$ is not empty, for $j = 1, \ldots, k$. 

 Let $H = \mathbb{Z}_{p_j}^{n_j} = \{h_1, \ldots, h_{r_j}\}$ be a fixed enumeration of elements of $H = \mathbb{Z}_{p_j}^{n_j} \subset G$. We denote by $a_1,\ldots,a_{q_j}$ a set of representatives of the left lateral classes of $G/\mathbb{Z}_{p_j}^{n_j}$ then, for each $i = 1, \cdots, q_j$, \,  $a_i h_1^{-1} = g_{l_1}, \ldots, a_i h_{r_j}^{-1} = g_{l_{r_j}}$ are elements of $G$.  Thus, if $y = (y_1,\ldots, y_r) \in A\left(\tilde{f}, \mathbb{Z}_{p_j}^{n_j} , G \right) $  then $$\tilde{f} (g_{l_1}\cdot(y_1,\ldots, y_r)) = \cdots = \tilde{f} (g_{l_{q_j}}\cdot(y_1,\ldots, y_r))$$ that is $$f(y_{\phi_{l_1(1)}}) = \cdots = f(y_{\phi_{l_{r_j}(1)}}).$$

Therefore, for each $j = 1, \ldots, k$, among the sides of ${\Delta}_N$ there are $r=q_{j}r_j$ mutually disjoint closed sides $\{{\sigma}_{i_1}$, ..., ${\sigma}_{i_{r_j}}\}_{i=1}^{q_{j}}$, such that 
$$f({\sigma}_{i_1})\cap \cdots \cap f({\sigma}_{i_{r_j}}) \neq \emptyset, $$ for each $i=1,\ldots, q_{j}$. \cqd

\vspace{1cm}

 Let us observe that since the $d$-dimensional Euclidean space $\R^d$ is homeomorphic to the interior of the closed $d$-dimensional ball, Theorem \ref{TverbergThm} hold also for maps into $\R^d$, and  we have the following {\it weak version of the topological Tverberg conjecture} or {\it topological Tverberg type theorem for any interger $r$.}
 
\begin{theorem}[Topological Tverberg type theorem for any interger $r$] \label{tv}Let $r\geq 2$, $d\geq 1$ be integers and $N=(r-1)(d + 1).$ Consider  $r = r_{1}\ldots r_k$ the prime factorization of $r$ and denote  $q_{j}=r/r_{j}$, $j=1,\ldots, k$. Then for any continuous map  $f:  {\Delta}_N \to \mathbb{R}^{d}$, for each $j = 1, \ldots, k$, there are $r=q_{j}r_{j}$ pairwise disjoint faces  $\{{\sigma}_{i_1}$, ..., ${\sigma}_{i_{r_j}}\}_{i=1}^{q_{j}}$  such that 
%${\sigma}_{11}$, ..., ${\sigma}_{1{r_i}}$;\ldots;${\sigma}_{j1}$, ..., ${\sigma}_{j{r_i}}$;\ldots; ${\sigma}_{{q_{i}}1}, \ldots, {\sigma}_{{q_{i}}{r_{i}}}$, such that 
$$f({\sigma}_{i_1})\cap \cdots \cap f({\sigma}_{i_{r_j}}) \neq \emptyset, \,\, \text{for each}\,\, i=1,\ldots, q_{j}.$$
\end{theorem} 

\begin{rema} \rm Let us note that if we consider $r$ a prime power  in Theorem \ref{tv}, we obtain the topological Tverberg theorem for prime powers.
\end{rema}

Now, by Theorem \ref{tv} and using similar method as in \cite{BFZ}, we have the following  {\it Generalized Van Kampen-Flores type  theorem for any interger $r$} or {\it a weak version of the Generalized Van Kampen-Flores   theorem}. In \cite[Theorem 4.2]{BFZ}, Blagojevic, Frick and Ziegler proved that the Generalized Van Kampen-Flores   theorem does not hold in general.
 
\begin{theorem}[Generalized Van Kampen-Flores type  theorem for any  $r$]\label{FloresThm}
Let $d\geq 1$ a natural number. Consider a natural number $r$ with prime factorization $r = r_{1}\cdots r_{k}$, $r_{1}<\cdots< r_{k}$, set $N = (r-1)(d+2)$  and let $l\geq [\frac{r-1}{r_{k}}d]$. Let $f:  \Delta_N \to \R^d$ be a continuous mapping. Then,  there are $r=q_{k}r_k$ pairwise disjoint faces  $\{{\sigma}_{i_1}, ..., {\sigma}_{i_{r_j}}\}_{i=1}^{q_{k}}$ of the $l$-th skeleton $\Delta^{(l)}_N$, such that 
$$f({\sigma}_{i_1})\cap \cdots \cap f({\sigma}_{i_{r_j}}) \neq \emptyset, \,\, \mbox{for each}\,\, i=1,\ldots, q_{k}.$$
\end{theorem}

\noindent {\it Proof}. Let $g: \Delta_N \to \R^{d+1}$ be a continuous function defined by $g(x) = (f(x), \dist(x, \Delta^{(l)}_N))$.  Then, we can apply  Theorem \ref{TverbergThm} to function $g$ which results in a collection of points $$x_{1_1}, ..., x_{1_{r_k}}; \ldots ; x_{i_1}, ..., x_{i_{r_k}}; \ldots; x_{{q_{k}}_1}, ..., x_{{q_{k}}_{r_{k}}},$$ 
 such that $\{x_{i_1}, ..., x_{i_{r_k}}\}_{i=1}^{q_{k}}$ are points in the   pairwise disjoint faces $\{{\sigma}_{i_1}, ..., {\sigma}_{i_{r_j}}\}_{i=1}^{q_{k}}$ with $f(x_{i_1}) = \cdots = f(x_{{i}_{r_{k}}})$ and $\dist(x_{i_1}, \Delta^{(l)}_N) = \cdots = \dist(x_{{i}_{r_{k}}}, \Delta^{(l)}_N) $, for each  $i=1,\ldots, q_{k}$.
We can suppose that the all ${\sigma}_{i_s}$'s are inclusion-minimal with the property that $x_{i_s} \in {\sigma}_{i_s}$, that is, ${\sigma}_{i_s}$ is the unique
face with $x_{i_s}$ in its relative interior. 

 Now, for each  $i=1,\ldots, q_{k}$ fixed, suppose that one of the faces ${\sigma}_{i_1},\ldots, {\sigma}_{i_{r_k}}$ is in $\Delta_{N}^{(l)}$, e.g. ${\sigma}_{i_1}$,  then $\dist(x_{i_{1}}, \Delta_{N}^{(l)})=0$, which imply that $\dist(x_{i_1}, \Delta^{(l)}_N) = \cdots = \dist(x_{{i}_{r_{k}}}, \Delta^{(l)}_N)=0$, and consequently, all faces ${\sigma}_{i_1},\ldots {\sigma}_{i_{r_k}}$ are in $\Delta_{N}^{(l)}$.
 
 Let us suppose the contrary, that no $\sigma_{i_{s}}$ is in $\Delta_{N}^{(l)}$, i.e., $\dim {\sigma}_{i_1}\geq l+1,\ldots, \dim{\sigma}_{i_{r_k}}\geq l+1$.  Since the faces ${\sigma}_{i_1},\ldots, {\sigma}_{i_{r_k}}$ are pairwise disjoint we have
  \begin{eqnarray*}
 N+1=| \Delta_{N}|&\geq& |{\sigma}_{i_1}|+\cdots +|{\sigma}_{i_{r_k}}| \\
 &\geq&r_{k}(l+2)\geq r_{k}(\left[\frac{r-1}{r_{k}}d\right]+2)\geq (r-1)(d+2)+2=N+2,\end{eqnarray*}
which is a contradiction and thus one of the faces ${\sigma}_{i_1},\ldots, {\sigma}_{i_{r_k}}$ is in $\Delta_{N}^{(l)}$  and consequently all faces ${\sigma}_{i_1},\ldots, {\sigma}_{i_{r_k}}$ are in $\Delta_{N}^{(l)}$.
\cqd

\vspace{1cm}

\begin{rema}\rm  Let us observe that if we consider $r$ a prime power  in Theorem \ref{FloresThm}, we obtain the Generalized Van Kampen-Flores  theorem for prime powers proved by Blagojevic, Frick and Ziegler in \cite[Theorem 3.2]{BFZ}.
\end{rema}

\vspace{0.5cm}

\small \noindent Denise de Mattos

\noindent {\it E-mail address} deniseml@icmc.usp.br

\noindent Universidade de S\~ao Paulo-USP-ICMC, Departamento de
Matemática, Caixa Postal 668, 13560-970, S\~ao Carlos-SP, Brazil

\vspace{0.5cm}

\noindent Edivaldo L. dos Santos

\noindent {\it E-mail address} edivaldo@dm.ufscar.br

\noindent Departamento de Matematica, Universidade Federal de Sao Carlos, Centro de Ciencias Exatas e Tecnologia, CP 676, CEP 13565-905, Sao Carlos - SP, Brazil.

\vspace{0.5cm}

\noindent Taciana O. Souza

\noindent {\it E-mail address} tacioli@ufu.br

\noindent           Faculdade de Matemáticatica, Universidade Federal de Uberlândia, Campus Santa Mônica - Bloco 1F - Sala 1F120, Av. João Naves de Avila, 2121, Uberlândia, MG, CEP: 38.408-100, Brazil

\end{document}